\documentclass[g5paper,twoside,final,10pt]{article}

\usepackage{amsfonts}
\usepackage{amsmath}
\usepackage{amsthm}
\usepackage{ifthen}
\newboolean{draftversion}

\newcommand{\note}[2][]{\message{(#1)}\ifthenelse{\boolean{draftversion}}%
	{\noindent{\em[#2]}}{}}

\def \P{\mathbb{P}}

\def \E{\mathbb{E}}

\def \eop {\hbox{}\nobreak\hfill
\vrule width 2mm height 2mm depth 0mm
\par \goodbreak \smallskip}

\def \eop {\hbox{}\nobreak\hfill \vrule width 2.0mm height 1.8mm depth 0mm
\par \goodbreak \smallskip}

\numberwithin{equation}{section}

\theoremstyle{plain}
\newtheorem{definition}{Definition}[section]
\newtheorem{theo}{Theorem}[section]

\newtheorem{lem}{Lemma}[section]
\newtheorem{remark}{Remark}[section]

\def \P{\mathbb{P}}

\def \E{\mathbb{E}}

\def \eop {\hbox{}\nobreak\hfill
\vrule width 2mm height 2mm depth 0mm
\par \goodbreak \smallskip}

\def \l{\left}
\def \r{\right}
\def \d{\mathrm{d}}

\def \P{\mathbb{P}}

\def \E{\mathbb{E}}

\def \eop {\hbox{}\nobreak\hfill
\vrule width 2mm height 2mm depth 0mm
\par \goodbreak \smallskip}

\def \eop {\hbox{}\nobreak\hfill \vrule width 2.0mm height 1.8mm depth 0mm
\par \goodbreak \smallskip}
\numberwithin{equation}{section}

\def \P{\mathbb{P}}

\def \E{\mathbb{E}}

\def \eop {\hbox{}\nobreak\hfill
\vrule width 2mm height 2mm depth 0mm
\par \goodbreak \smallskip}

\def \eop {\hbox{}\nobreak\hfill
\vrule width 1.4mm height 1.4mm depth 0mm
\par \goodbreak \smallskip}

\def \P{\mathbb{P}}

\def \E{\mathbb{E}}

\def \eop {\hbox{}\nobreak\hfill\vrule width 2mm height 2mm depth 0mm
\par \goodbreak \smallskip}

\def \eop {\hbox{}\nobreak\hfill \vrule width 2.0mm height 1.8mm depth 0mm
\par \goodbreak \smallskip}

\title{A mixed relaxed singular maximum principle for linear SDEs with random coefficients}
\begin{document}
\author{Daniel Andersson\thanks{Department of
Mathematics, Royal Institute of Technology, S-100 44 Stockholm,
Sweden. danieand@math.kth.se}}
\maketitle

\begin{abstract} \noindent We study singular stochastic control of a two dimensional stochastic differential equation, where the first component is linear with random
and unbounded coefficients. We derive existence of an optimal relaxed control and necessary conditions for optimality in the form of a mixed relaxed-singular maximum principle in a global form.
A motivating example is given in the form of an optimal investment and consumption problem with transaction costs, where we consider a portfolio with a continuum of bonds and where the
 portfolio weights are modeled as measure-valued processes on the set of times to maturity.
\bigskip
\bigskip


\end{abstract}


\section{Introduction}
The objective of this paper is to derive necessary conditions for optimality in mixed relaxed-singular stochastic control problems. That is, the control has two parts: one absolutely continuous and one singular. The relaxation is performed by replacing the absolutely continuous control with a control that takes values on the set of probability measures.
The state process is a solution to a two dimensional stochastic differential equation (SDE). The first component is a linear SDE whose coefficients are random and not necessarily bounded. The second component is a general non-linear SDE whose coefficients have bounded derivatives.

\medskip\noindent A motivating example is the following optimal consumption-investment problem. We consider a market with two investment opportunities, a stock and a portfolio of bonds. The bonds are non-defaultable, i.e. financial contracts that are bought today and pay a fixed amount at some
 future time, called the maturity time. At each time $t$, the investor is allowed to buy bonds with any time to maturity in $U$, where $U$ is a subset of $\mathbb{R}_+$. The relative portfolio weights are therefore modeled as a probability measure on $U$, reflecting the proportion invested in bonds with different maturities. Modeling the prices of the bonds as SDEs, we may write down the value $x_t$ of the bondportfolio as an SDE of the form (see Section \ref{example} below)
\begin{align*}
x_t=x_0+\int_0^tx_s\int_U\l(r_s^0-v_s\left(u
\right)\Theta_s\right)q_s\left(\d u
\right)\d s+\int_0^tx_s\int_Uv_s\left(u
\right)q_s\left(\d u
\right)\d B^x_s,
\end{align*}
where $x_0$ is the initial capital, $B^x_t$ is a Brownian motion and $q_t$ is a probability measure on $U$. Further, $r_t^0$ is
the short rate, $v_t$ is the integrated volatility process of the bond prices and $\Theta_t$ is the so called market price of risk.

\medskip\noindent The price of a share of a stock is modeled as a geometric Brownian motion and the value of the investment in the stock at time $t$ is then given by
\begin{align*}
y_t=y_0+\int_0^t\lambda y_s \d s+ \int_0^t\rho y_s \d B^y_s,
\end{align*}
where $y_0$ is the initial capital, $\lambda$ and $\rho$ are constants and $B^y_t$ is a Brownian motion independent of $B^x_t$. Note that the prices of bonds and the stock are modeled under the physical measure and not under the so-called risk neutral probability measure.

\medskip\noindent Denoting throughout the paper
\begin{align*}
\varphi_t\l(q_t\r)=\int_U \varphi_t\l(u\r)q_t\l(\d u\r)
\end{align*}
for any function $\l(\omega,t,u\r)\mapsto \varphi_t(u)$, the position at time $t$ for an investor with these two investment possibilities is $\l(x_t,y_t\r)$ given by
 \begin{align*}
 x_t=&x_0+\int_0^tx_s\left(r_s^0-v_s\left(q_s
\right)\Theta_s-c_s\right)\d s+\int_0^tx_sv_s\left(q_s
\right)\d B^x_s+\left(1-K_1\right) \xi^x_t-\xi^y_t,\\
 y_t=&y_0+\int_0^t\lambda y_s\d s+\int_0^t\rho y_s\d B^y_s-\xi^x_t+\left(1-K_2\right)\xi^y_t,
 \end{align*}
where the consumption process $c_t$ is required to take values in some compact subset of $\mathbb{R}_+$ and $\xi_t=\l(\xi_t^x,\xi_t^y\r)$ is nondecreasing and left continuous with right limits. The value of the bonds sold to buy stocks is recorded by $\xi_t^y$, and $\xi_t^x$ records the value of the stocks sold to buy bonds. 
The constants $0\leq K_1,K_2<1$ account for the proportional transaction costs incurred whenever money is moved between the stock and the bonds. The position in the stock, $y_t$, is independent of the absolutely continuous control, reflecting a buy-and-hold strategy, while the strategy in the bond market involves continuous rebalancing of the portfolio.

\medskip\noindent The objective of the investor is then to choose a consumption/investment strategy to minimize some cost functional
\begin{align}
J\l(\mu,\c,\xi\r)=\E\l(\int_0^Th\l(t,x_t,y_t,\mu_t,c_t\r)\d t+\int_0^Tk_t\d\xi_t+g\l(x_T,y_T\r)\r),\label{costintr}
\end{align}
That is, the objective is to optimally choose three adapted processes $\l(c_t,q_t,\xi_t\r)$ such that (\ref{costintr}) is minimized. 

\medskip\noindent This is an example of a singular stochastic control problem with some non standard characteristics. Firstly, the state process is a linear SDE with random coefficients, and where $r_t^0$ and $\Theta_t$ cannot in general
be assumed to be bounded. Secondly, the absolutely continuous part of the control is extended from the action space $U$ to the space $\mathcal{P}(U)$ of probability measures on $U$. 

\medskip\noindent This motivates us to study control problems of the form
\begin{align*}
x_t=&x_0+\int_0^tb^x\left(s,x_s,\mu_t\right)\d s+\int_0^t\sigma^x\left(s,x_s,\mu_t\right)\d B_s+\int_0^tG^x_s\d \xi_s,\\
y_t=&y_0+\int_0^tb^y\l(s,y_s\r)\d s+\int_0^t\sigma^y\l(s,y_s\r)\d B_s+\int_0^tG^y_s\d \xi_s.
\end{align*}
The control is a process $\mu_t$ taking values in the space of probability measures on the action space $U$ and a nondecreasing process $\xi_t$, left continuous with limits on the right. $B_t$ is a $d$-dimensional Brownian motion,
 $x_0$ denotes the
 initial state, $b^x$ and $\sigma^x$ are random coefficients of the form
\begin{align*}
b^x\left(t,x,u,\omega\right)=\upsilon_t\left(u,\omega\right)+\phi_t\left(u,\omega\right)x,\\
\sigma^x\left(t,x,u,\omega\right)=\chi_t\left(u,\omega\right)+\psi_t\left(u,\omega\right)x,
\end{align*}
for given stochastic processes $\upsilon$, $\phi$, $\chi$ and $\psi$ taking values in the space of continuous functions on $U$. Further, $b^y$ and $\sigma^y$ are deterministic functions and 
the cost functional, which is to be minimized, is of the form
\begin{align*}
 J\left(\mu_t,\xi_t\right)=\E\left(\int_0^Th\left(t,x_t,y_t,\mu_t\right)\d t+\int_0^Tk_t\d\xi_t+g\left(x_T,y_T\right)\right).
 \end{align*}
\noindent This paper contains two main results. The first one, Theorem \ref{existence}, establishes existence of an optimal relaxed control which is derived using a similar scheme as in \cite{anderssondjehiche}. The main tools in the proof are tightness and Skorohod's selection theorem. The second main result, Theorem \ref{maxprinciple}, suggests necessary conditions for optimality that are given in form of a relaxed maximum principle. We follow the scheme in 
\cite{bahlalimezerdi}, where a stochastic maximum principle of second order type (i.e.~two adjoint processes) is obtained by performing a spike perturbation on the absolutely continuous control. However, by only considering relaxed controls, i.e.~extending the control to the space of probability measures on the action space, $\mathcal{P}(U)$, which is a convexification of the action space, it allows us to use a convex perturbation although the action space is not convex. Ultimately this leads to a maximum principle in global form, i.e.~just one adjoint process. 

\medskip\noindent Under the usual assumptions on the coefficients in the SDE, i.e.~deterministic functions of $\left(t,x,u\right)$, Lipschitz continuous and with linear growth in $x$, a maximum principle for stochastic (strict) control problems where the
control enters the diffusion coefficient and the control set is not convex was established in \cite{peng}. Since the action space is not convex, a spike perturbation method was applied which led to a second order maximum principle. 

\medskip\noindent As for singular stochastic control problems, they have been studied by many authors, see \cite{haussmannsuo} and the references therein. These papers mainly focus on the dynamic programming principle. The first stochastic maximum principle for singular control problems is obtained in \\ \cite{cadenillashaussmann}, where they assume linear dynamics with random but bounded coefficients, convex cost criterion and convex state constraint. In \cite{bahlalietal} the maximum principle is extended to include relaxed controls, under the assumption that the diffusion coefficient is independent of the control.

\medskip\noindent The paper is organized as follows. In Section \ref{formulation}, we formulate the mixed relaxed-singular control problem for our linear SDEs. In Section \ref{existenceresult} we prove existence of an optimal control, while in Section \ref{necessary}, necessary conditions for optimality are given in form of a relaxed maximum principle. In Section \ref{example}, we apply these results to formulate a maximum principle for the optimal investment/consumption problem.
\section{Formulation of the problem}\label{formulation}

Let $T>0$ be a fixed time horizon and $\left( \Omega ,%
\mathcal{F},\mathcal{F}_t,\P\right) $ be a filtered
probability space satisfying the usual conditions, on which a $d-$%
dimensional Brownian motion $\l\{B_t\r\}_{t\in[0,T]}$ is defined. We
assume that $\left( \mathcal{F}_{t}\right) _{t\in[0,T]}$ is the natural filtration
of $B_t$ augmented by $\P-$null sets of $\mathcal{F}.$

\medskip\noindent Consider the following sets, $U$ is a compact subset of $\mathbb{R}%
^{d} $ and $A=\left( \left[ 0,\infty \right) \right) ^{d}.$
Let $\mathcal{U}$ the class of measurable, adapted processes $u:\left[ 0,T\right]
\times \Omega \longrightarrow U$ and $\mathcal{A}$ the class of measurable,
adapted processes $\xi :\left[ 0,T\right] \times \Omega \longrightarrow
A$ such that $\xi $ is of bounded variation, nondecreasing
left continuous with right limits, $\xi _{0}=0$ and $\E\l|\xi_T\r|^p<\infty$ for any $p\geq 1$.

\medskip\noindent We define a two dimensional controlled SDE on $\left(\Omega,\mathcal{F},\mathcal{F}_t,\P\right)$, with absolutely continuous control $u_t$ and singular control $\xi_t$:
\begin{subequations}\label{strictstateequation}
\begin{align}
x_t=&x_0+\int_0^tb^x\left(s,x_s,u_s\right)\d s+\int_0^t\sigma^x\left(s,x_s,u_s\right)\d B_s+\int_0^tG^x_s\d \xi_s,\\
y_t=&y_0+\int_0^tb^y\l(s,y_s\r)\d s+\int_0^t\sigma^y\l(s,y_s\r)\d B_s+\int_0^tG^y_s\d \xi_s.
\end{align}
\end{subequations}
where $\l(x_0,y_0\r)\in\mathbb{R}^2$ is the initial state, $G_t^x,G^y_
t:[0,T]\mapsto\mathbb{R}^d$, $b^y:[0,T]\times\mathbb{R}\mapsto\mathbb{R}$ and $\sigma^y:[0,T]\times\mathbb{R}\mapsto\mathbb{R}^d$.
 Furthermore, the coefficients $b^x$ and $\sigma^x$ are given by
\begin{subequations} \label{linearcoeff}
\begin{align}
b^x\left(t,x,u,\omega\right)=\upsilon_t\left(u,\omega\right)+\phi_t\left(u,\omega\right)x  \\
\sigma^x\left(t,x,u,\omega\right)=\chi_t\left(u,\omega\right)+\psi_t\left(u,\omega\right)x,
\end{align}
\end{subequations}
where $\upsilon:[0,T]\times U\times\Omega\mapsto \mathbb{R}$, $\phi: [0,T]\times U\times\Omega
\mapsto \mathbb{R}$, $\chi: [0,T]\times U\times\Omega\mapsto \mathbb{R}^d$ and
$\psi:[0,T]\times U\times\Omega\mapsto \mathbb{R}^d$ are $\mathcal{F}_t$-adapted
processes. 

\medskip\noindent The cost functional is given by
\begin{align}
J\left(u,\xi\right)=\E\left(\int_0^Th\left(t,x_t,y_t,u_t\right)\d t+\int_0^tk_t\d\xi_t+g\left(x_T,y_T\right)\right), \label{costfunctional1}
\end{align}
and the objective is to minimize $J$ over the set of admissible controls. A control $\l(u^*,\xi^*\r)$ is called optimal if it satisfies $J\left(u^*,\xi^*\right)=\inf\{J\left(u,\xi\right);u\in\mathcal{U},\xi\in\mathcal{A}\}$. If also $\l(u^*,\xi^*\r)\in
\mathcal{U}\times\mathcal{A}$, it is called a strict optimal control.

\medskip\noindent This kind of control problems is often formulated in the so-called relaxed form, due to the fact that a strict optimal control may fail to exist (see e.g.~\\\cite{boualemetal} and \cite{lou}). Instead one embeds the strict
controls in a wider class of controls that takes values in probability measures on $U$ rather than on $U$ itself. Also, a solution to a relaxed control problem is a weak one, i.e.~the probability space, equipped with the a priori given stochastic processes, is part of the solution.

\medskip\noindent Let $\mathcal{P}( U)$ be the space of probability measures on $U$. If
$\mu_t(\d u)$ is a stochastic process taking values in $\mathcal{P}(U)$, we denote by $\mathcal{L}([0,T],U)$ the space of the (Radon) measure-valued processes
$\d \lambda_t(u)=\mu_t(\d u)\d t$. If a probability space $\left(\Omega, \mathcal{F}, \P\right)$ is given, then we denote $\boldsymbol{M}(\Omega)$ the space of all $\mathcal{F}_t$-adapted processes $\mu_t(\d u)$ taking values in $\mathcal{P}( U)$. Further, we denote by $\boldsymbol{L}(\Omega)$ the space of all $\mathcal{L}([0,T],U)$-valued $\mathcal{F}_t$-adapted processes. It can be shown that there is a one-to-one correspondence between  $\boldsymbol{M}(\Omega)$ and $\boldsymbol{L}(\Omega)$, and that $\mathcal{L}([0,T],U)$ is a compact metric space. For further discussion, see \cite{mayong}. 

\medskip\noindent By expanding the set of controls from $\mathcal{U}$ to $\boldsymbol{M}$,
  the state equation is defined as
\begin{subequations}\label{relaxedstateeq}
\begin{align}
x_t=&x_0+\int_0^tb^x\left(s,x_s,\mu_s\right)\d s+\int_0^t\sigma^x\left(s,x_s,\mu_s\right)\d B_s+\int_0^tG^x_s\d \xi_s,\\
y_t=&y_0+\int_0^tb^y\l(s,y_s\r)\d s+\int_0^t\sigma^y\l(s,y_s\r)\d B_s+\int_0^tG^y_s\d \xi_s.
\end{align}
\end{subequations}

\medskip\noindent We make the following assumptions regarding the state equation and the cost functional.
\begin{itemize}
\item[(A.1)] $\varphi_t\left(u,\omega\right)$ is continuous in $\left(t,u\right)$, where $\varphi$ stands for one of the processes $\upsilon, \phi, \chi, \psi$.
\item[(A.2)] $\psi$ is uniformly bounded in $\mathbb{R}^+\times U\times\Omega$.
\item[(A.3)] For any $p\in(-\infty,\infty)$, it holds that
\begin{align*}
\underset{u\in U}{\sup}\E\left(\exp\left(p\int_0^T\phi_t\left(u\right)\d t\right)\right)<\infty.
\end{align*}
\item[(A.4)] For any $p\geq 1$, it holds that
\begin{align*}
\underset{u\in U}{\sup}\E\left(\int_0^T|\upsilon_t\left(u\right)|^p\d t\right)<\infty\,\,\, \mbox{and}\,\,\, \underset{u\in U}{\sup}\E\left(\int_0^T|\chi_t\left(u\right)|^p\d t\right)<\infty.
\end{align*}
\item[(A.5)] The functions $b^y,\sigma^y$ are continuously differentiable in $y$ with bounded derivatives. The functions $G^x,G^y$ are positive, continuous and bounded.
\item[(A.6)] The functions $g$ and $h$ are continuously differentiable in $\l(x,y\r)$. The function $g$ and its derivative are bounded and Lipschitz continuous in $\l(x,y\r)$. The function $h$ and its derivative are bounded, continuous in $u$ and Lipschitz continuous in $\l(x,y\r)$. The function $k$ is continuous and bounded.
 \end{itemize}

\begin{remark}
Assumption $(A.6)$ is made here out of convenience. It is possible to consider $h$ and $g$ that are only Lipschitz continuous in $\l(x,y\r)$, without any assumption about differentiability. See \cite{kbahlalietal}. 
\end{remark}

\begin{definition} \label{defrelaxedcontrol}
A relaxed control is the term $\mathcal{A}=\left(\Omega, \mathcal{F}, \mathcal{F}_t, \P, \boldsymbol{B}_t, \mu_t, x_t, y_t, \xi_t \right) $, where
\begin{itemize}
\item[$(i)$] $\left(\Omega, \mathcal{F}, \mathcal{F}_t, \P\right)$  is a filtered probability
space;
\item[$(ii)$] $\boldsymbol{B}_t= \left(B_t,\upsilon_t, \phi_t, \chi_t, \psi_t\right)$, in which $\l\{B_t\r\}$ is a $d$-dimensional $\mathcal{F}_t$-Brownian motion and
$\upsilon_t, \phi_t, \chi_t, \psi_t$ are $\mathcal{F}_t$-adapted stochastic processes satisfying \textnormal{(A.1)-(A.4)};
\item[$(iii)$] $\mu_t\in
\boldsymbol{M}(\Omega)$, $\xi_t\in \mathcal{A}$;
\item[$(iv)$] $\l(x_t,y_t\r)$ is $\mathcal{F}_t$-adapted and satisfies \textnormal{(\ref{relaxedstateeq})}.
\end{itemize}
\end{definition}

\medskip
\noindent We denote by $\mathcal{R}$ the set of all relaxed controls. The cost functional corresponding to the control $\mathcal{A}$
 is defined as
\begin{align}
J\left(\mathcal{A}\right)=\E\left(\int_0^Th\left(t,x_t,y_t,\mu_t\right)\d t+\int_0^tk_t\d\xi_t+g\left(x_T,y_T\right)\right), \label{costfunctional}
\end{align}
and a relaxed control $\mathcal{A}^*$ is optimal if $J\left(\mathcal{A}^*\right)=\inf\{J\left(\mathcal{A}\right);\mathcal{A}\in\mathcal{R}\}$. It is well known that $\mathcal{U}$ may be embedded into
$\mathcal{R}$, since any strict ($U$-valued) control process $u_t$ can be represented as a relaxed control by setting $\mu_t(\d u)=\delta_{u_t}(\d u)$.

\medskip
\noindent Throughout the rest  of the paper we will not specify that properties hold $\P$-a.s. when it is clear from the context. Further, we denote for any process $\varphi_t$, \begin{align*}
\vert\varphi\vert^{*,p}_T=\sup_{t\in[0,T]}|\varphi_t|^p.
\end{align*}
\section{Existence of an optimal control}\label{existenceresult}

In this section we shall establish the existence of an optimal relaxed control. To achieve this, we consider a minimizing sequence of controls $\mathcal{A}^{(k)}\in\mathcal{R}$ for the cost functional $J$, i.e.
\begin{align*}
\inf\{J(\mathcal{A}),\mathcal{A}\in\mathcal{R}\}=\lim_{k\to\infty}J(\mathcal{A}^{(k)}),
\end{align*}
and show that a limit $\mathcal{A}$ exists and fulfills $(i)-(iv)$ in Definition \ref{defrelaxedcontrol}. The main tools are tightness of the processes and Skorohod's
 Selection Theorem.

\medskip\noindent Given a relaxed control $\mathcal{A}=(\Omega,\mathcal{F},\mathcal{F}_t,\P,\boldsymbol{B}_t,\mu_t,x_t,y_t,\xi_t)$, there exists a unique strong solution to the equation given by
(\ref{linearcoeff}) and (\ref{relaxedstateeq}).  Moreover, with the assumptions (A.1)-(A.5) we can prove by standard methods that
$\l(x_t,y_t\r)$ has the following properties:
For any $p\geq1$ we have
\begin{align}
\E\l|x_\cdot\r|_T^{*,p} + \E\l|y_\cdot\r|_T^{*,p} < \infty. \label{momentsde}
\end{align}

\noindent We endow the space $\mathcal{A}$ with the pseudopath topology, cf.~\cite{haussmannsuo}.
The pseudopath identifies two functions if and only if they are equal (Lebesgue) almost everywhere. Under the pseudopath topology, $\mathcal{A}$ is a seperable metric space and convergence in the pseudopath topology is just convergence in measure, i.e.~for $\xi^n,\xi\in\mathcal{A}$, $\xi^n\to \xi$ if and only if
\begin{align}
\int_0^Tf(t)~\d \xi_t^n\to \int_0^Tf(t)~\d \xi_t \label{convxi}
\end{align}
for any $f\in C\l([0,T],\mathbb{R}^d\r)$.


\medskip
\noindent Let $\mathcal{A}^{(k)}=(\Omega^{(k)},\mathcal{F}^{(k)},\P^{(k)},\mathcal{F}_t^{(k)},\boldsymbol{B}_t^{(k)},\mu_t^{(k)},x_t
^{(k)}, y_t^{(k)}, \xi_t^{(k)})$ be a minimizing sequence, i.e.
\begin{align*}
J(\mathcal{A}^{(k)})\to
\inf\{J(\mathcal{A});\mathcal{A}\in\mathcal{U}^R\}, ~ \textrm{as} ~ k\to \infty. \end{align*}
We will derive tightness for this sequence. First, we introduce the notation
\begin{align*}
\lambda_t^{(k)}(A)&=\int_0^t\int_A\mu_s^{(k)}(du)ds,
\end{align*}
for any Borel set $A\subset U$, and define the continuous processes
\begin{align*}
m_t^{(k)}(x)&=x^{(k)}_t-\int_0^tG^x_s\d \xi_s^{(k)},\\
m_t^{(k)}(y)&=y^{(k)}_t-\int_0^tG^y_s\d \xi_s^{(k)}.
\end{align*}

\begin{lem}\label{tight}
The sequence
$(\boldsymbol{B}_t^{(k)},\lambda_t^{(k)},m(x)_t^{(k)}, m(y)^{(k)}_t, \xi_t^{(k)})$ is tight in\\
$\left(C([0,T])\times C([0,T]\times U)^4\right)\times\mathcal{L}([0,T]\times U)\times C([0,T])\times C([0,T])\times \mathcal{A}$.
 \end{lem}
\noindent\textbf{Proof.}
$\boldsymbol{B}_t^{(k)}$ is tight since the processes induce the same measure for every $k$. Also,
$\lambda_t^{(k)}$ is tight because $\mathcal{L}([0,T]\times U)$ is compact. Moreover, by (A.1)-(A.5), it is readily seen that there exists a constants $K_1, K_2>0$ such that
\begin{align*}
\E^{(k)}(m(x)_t^{(k)}-m(x)_s^{(k)})^4\leq K_1|s-t|^2,\\
\E^{(k)}(m(y)_t^{(k)}-m(y)_s^{(k)})^4\leq K_2|s-t|^2,
\end{align*}
for all $t,s \in [0,T]$, for all $k$, where $\E^{(k)}$ is the expectation under $\P^{(k)}$.
 Hence the Kolmogorov condition is fulfilled (see e.g. \cite{yongzhou}, Theorem 2.14.) and $\l(m(x)_t^{(k)},m(y)_t^{(k)}\r)$ is tight.

\medskip\noindent As for the tightness of $\xi_t^{(k)}$, we proceed as in \cite{haussmannsuo} and define the set
\begin{align*}
V_M=\l\{\xi\in \mathcal{A}^k : \l|\xi_T\r|\leq M\r\}.
\end{align*}
$V_M$ is then compact for any constant $M>0$. Further, define
\begin{align*}
\mathcal{R}^\lambda=\l\{\mathcal{A}\in\mathcal{R} : J(\mathcal{A})\leq\lambda\r\},
\end{align*}
where $\lambda$ is chosen so that $\mathcal{R}^\lambda$ is nonempty. Obviously we can restrict the minimizing sequence to $\mathcal{R}^\lambda$. It also holds that
\begin{align*}
\underset{M\to\infty}{\lim}~\underset{\mathcal{A}\in \mathcal{R}^\lambda}{\inf}~\P\l(\l|\xi_T\r|\leq M\r)=1,
\end{align*}
see \cite{haussmannsuo}, Proposition 3.4.
Thus, for any given $\varepsilon>0$ there exists a compact set $V_M$ such that for all $\P^{(k)}, \xi_t^{(k)}\in \mathcal{R}^\lambda$
\begin{align*}
\P^{(k)}\l(\xi_\cdot^{(k)} \in V_M\r)\geq 1-\varepsilon.
\end{align*}
Thus, the sequence $\xi_t^{(k)}$ is tight.
\eop

\begin{theo}\label{existence}
Under the assumptions $(A.1)-(A.6)$, the relaxed control problem admits an optimal solution.
\end{theo}
\noindent\textbf{Proof.}
 By Lemma \ref{tight}, the sequence $(\boldsymbol{B}_t^{(k)},\lambda_t^{(k)},m(x)_t^{(k)}, m(y)^{(k)}_t, \xi_t^{(k)})$ is tight and thus by Skorohod's selection theorem, there exists a probability space
$(\hat{\Omega},\hat{\mathcal{F}},\hat{\P})$, on which is defined a sequence of processes $(\hat{\boldsymbol{B}}_t^{(k)},\hat{\lambda}_t^{(k)},m(\hat{x})_t^{(k)}, m(\hat{y})_t^{(k)}, \hat{\xi}_t^{(k)})$
identical in law to $(\boldsymbol{B}_t^{(k)},\lambda_t^{(k)},m(x)_t^{(k)},m(y)_t^{(k)},\xi_t^{(k)})$ and converging $\hat{\P}$-a.s. to $(\hat{\boldsymbol{B}_t},\hat{\lambda}_t,m(\hat{x})_t, m(\hat{y})_t,\hat{\xi}_t)$.

\medskip\noindent Moreover, as in \cite{anderssondjehiche}, Lemma 3.1 we may assume that $\mu_t^{(k)}(A)$ has continuous paths for each Borel set $A$. Thus, by Lemma
2.1 in \cite{mayong} the processes $\hat{\mu}_t^{(k)}$ corresponding to $\hat{\lambda}_t^{(k)}$ have the same law as $\mu_t^{(k)}$ since $\mu_t^{(k)}$ has continuous paths.

\medskip\noindent We have then shown that there exists a limit $\hat{\mathcal{A}}$ of the minimizing
 sequence $\hat{\mathcal{A}^{(k)}}$ which satisfies $(i)-(iii)$ in Definition \ref{defrelaxedcontrol}. It remains to show that $(iv)$ also holds.
To this end let $\hat{b^x}$ and $\hat{\sigma^x}$ be the processes defined by \textnormal{(\ref{linearcoeff})}, then we can prove as in \cite{mayong} Lemma 3.3, that
\begin{align*} \hat{b^x}(t,\hat{x}_t,\hat{\mu}_t^{(k)})\stackrel{w}\longrightarrow \hat{b^x}(t,\hat{x}_t,\hat{\mu}_t),\\
  \hat{\sigma^x}(t,\hat{x}_t,\hat{\mu}_t^{(k)})\stackrel{w}\longrightarrow \hat{\sigma^x}(t,\hat{x}_t,\hat{\mu}_t),
 \end{align*}  in $L^2([0,T]\times\Omega)$, as $ k\to
\infty$.
 Using this result and the fact that \\$(\hat{\boldsymbol{B}}_t^{(k)},\hat{\lambda}_t^{(k)},m(\hat{x})_t^{(k)}, m(\hat{y})_t^{(k)}, \hat{\xi}_t^{(k)})\to(\hat{\boldsymbol{B}}_t,\hat{\lambda}_t,m(\hat{x})_t, m(\hat{y})_t,\hat{\xi}_t)$~$\hat{\P}$-a.s., as $k\to\infty$, we can prove that the limit ~$\l(\hat{x}_t,\hat{y}_t\r)$ of $\l(\hat{x}_t^{(k)},\hat{y}_t^{(k)}\r)$ satisfies
\begin{align}
\hat{x}_t&=x_0+\int_0^t\hat{b^x}\l(s,\hat{x}_s,\hat{\mu}_s\r)\d s+\int_0^t\hat{\sigma^x}\l(s,\hat{x}_s,\hat{\mu}_s\r)\d \hat{B_s}+\int_0^tG^x_s\d\hat{\xi}_s,\\
\hat{y}_t&=y_0+\int_0^tb^y\l(s,\hat{y}_s\r)\d s+\int_0^t \sigma^y\l(s,\hat{y}_s\r) \d \hat{B_s}+\int_0^tG^y_s\d\hat{\xi}_s, \label{sdey}
\end{align}
i.e. that
\begin{align*}
&\hat{\E}\l|\hat{x}_\cdot-x_0-\int_0^\cdot \hat{b^x}\l(s,\hat{x}_s,\hat{\mu}_s\r)\d s-\int_0^\cdot\hat{\sigma^x}\l(s,\hat{x}_s,\hat{\mu}_s\r)\d \hat{B_s}-\int_0^\cdot G^x_s\d\hat{\xi}_s\r|_T^{*,2}\\
+&\hat{\E}\l|\hat{y}_\cdot-y_0-\int_0^\cdot b^y\l(s,\hat{y}_s\r)\d s-\int_0^\cdot\sigma^y\l(s,\hat{y}_s\r) \d \hat{B_s}-\int_0^\cdot G^y_s\d\hat{\xi}_s\r|_T^{*,2}=0.
\end{align*}
As for $\hat{x}_t$, we note that by (\ref{convxi}) the sequence $m(\hat{x})_t^{(k)}=\hat{x}_t^{(k)}-\int_0^t G^x_s\d\hat{\xi}_s^{(k)}$ converges to $m(\hat{x})_t=\hat{x}_t-\int_0^t G^x_s\d\hat{\xi}_s$ $\hat{\P}$-a.s. in $C([0,T])$ and by uniform integrability
\begin{align*}\hat{\E}\l|\hat{x}_\cdot^{(k)}-\int_0^\cdot G^x_s\d\hat{\xi}_s^{(k)}-\hat{x}_\cdot+\int_0^\cdot G^x_s\d\hat{\xi}_s\r|^{*,2}_T\to 0.
\end{align*}
Further, by applying the same proof as in \cite{anderssondjehiche}, Theorem 3.1, we can prove that
\begin{align*}
&\hat{\E}\l|\int_0^\cdot\hat{b^x}^{(k)}\l(s,\hat{x}^{(k)}_s,\hat{\mu}^{(k)}_s\r)\d s-\int_0^\cdot\hat{b^x}\l(s,\hat{x}_s,\hat{\mu}_s\r)\d s\r|_T^{*,2}\to 0,~\textrm{and}\\
&\hat{\E}\l|\int_0^\cdot\hat{\sigma^x}^{(k)}\l(s,\hat{x}^{(k)}_s,\hat{\mu}^{(k)}_s\r)\d \hat{B}^{(k)}_s-\int_0^\cdot\hat{\sigma^x}\l(s,\hat{x}_s,\hat{\mu}_s\r)\d \hat{B_s}\r|_T^{*,2}\to 0,
\end{align*}
as $k\to \infty$.
 Now, since $\hat{x}_t^{(k)}$ satisfies
\begin{align*}
\hat{x}_t^{(k)}=x_0+\int_0^t\hat{b^x}^{(k)}\l(s,\hat{x}^{(k)}_s,\hat{\mu}^{(k)}_s\r)\d s+\int_0^t\hat{\sigma^x}^{(k)}\l(s,\hat{x}^{(k)}_s,\hat{\mu}^{(k)}_s\r)\d \hat{B}^{(k)}_s+\int_0^tG^x_s\d\hat{\xi}_s^{(k)},
\end{align*}
for every $k$, we may write
\begin{align*}
&\hat{\E}\l|\hat{x}_\cdot-x_0-\int_0^\cdot \hat{b^x}\l(s,\hat{x}_s,\hat{\mu}_s\r)\d s-\int_0^\cdot\hat{\sigma^x}\l(s,\hat{x}_s,\hat{\mu}_s\r)\d \hat{B}_s-\int_0^\cdot G^x_s\d\hat{\xi}_s\r|_T^{*,2}\\
\leq&~ \hat{\E}\l|\hat{x}_\cdot-G^x_s\d\hat{\xi}_s-\hat{x}_\cdot^{(k)}+\int_0^tG^x_s\d\hat{\xi}_s^{(k)}\r|^{*,2}_T\\
+&\hat{\E}\l|\int_0^t\hat{b^x}^{(k)}\l(s,\hat{x}^{(k)}_s,\hat{\mu}^{(k)}_s\r)\d s-\int_0^t\hat{b^x}\l(s,\hat{x}_s,\hat{\mu}_s\r)\d s\r|_T^{*,2}\\
+&\hat{\E}\l|\int_0^t\hat{\sigma^x}^{(k)}\l(s,\hat{x}^{(k)}_s,\hat{\mu}^{(k)}_s\r)\d \hat{B}^{(k)}_s-\int_0^t\hat{\sigma^x}\l(s,\hat{x}_s,\hat{\mu}_s\r)\d \hat{B}_s\r|_T^{*,2}.
\end{align*}
The conclusion follows by letting $k\to\infty$. Similarly (more easily, in fact) we can show that $\hat{y}_t$ satisfies (\ref{sdey}).

\medskip\noindent Finally, since $f$, $g$ and $k$ are continuous and bounded, it is readily seen that
\begin{align*}
&J\big(\hat{\mathcal{A}}^{(k)}\big)=\hat{\E}\l(\int_0^Tf\l(s,\hat{x}^{(k)}_s,\hat{y}^{(k)}_s,\hat{\mu}^{(k)}_s\r)\d s+\int_0^Tk_s\d \hat{\xi}^{(k)}_s\r)\\
&\to \hat{\E}\l(\int_0^Tf\l(s,\hat{x}_s,\hat{y}_s,\hat{\mu}_s\r)\d s+\int_0^Tk_s\d \hat{\xi}_s\r)=J\big(\hat{\mathcal{A}}\big).
\end{align*}
Since $\hat{\mathcal{A}}^{(k)}$ is a minimizing sequence , we conclude that
\begin{align*}
J\big(\hat{\mathcal{A}}\big)=\underset{\mathcal{A}\in\mathcal{R}}{\inf}J\l(\mathcal{A}\r).
\end{align*}
Therefore, $\hat{\mathcal{A}}=(\hat{\Omega},\hat{\mathcal{F}},\hat{\P},\hat{\mathcal{F}}_t,\hat{\boldsymbol{B}}_t,\hat{\mu}_t,\hat{x}_t,\hat{y}_t,\hat{\xi}_t)$ is an optimal relaxed control.
\eop


\section{A relaxed maximum principle} \label{necessary}

\subsection{Preliminary results}
Let $\left(\mu,\xi\right)\in \mathcal{R}$ be an optimal control and $x_t$ the corresponding state trajectory.
Since $\mathcal{P}(U)$ is a convexification of the action space $U$ we may introduce the following convex perturbation of $\left(\mu,\xi\right)$.
\begin{align*}
\mu_t^\theta=\mu_t+\theta\left(q_t-\mu_t\right),\\
\xi_t^\theta=\xi_t+\theta\left(\eta_t-\xi_t\right),
\end{align*}
for $0<\theta<1$ and $\l(q,\eta\r)\in\mathcal{R}$. Now, since $\mathcal{P}(U)$ is a convex space, it is readily seen that $\l(\mu^\theta,\xi^\theta\r)\in\mathcal{R}$. 
Further, we define
\begin{align*}
J_1=J\left(\mu^\theta,\xi^\theta\right)-J\left(\mu^\theta,\xi\right),\\
J_2=J\left(\mu^\theta,\xi\right)-J\left(\mu,\xi\right).
\end{align*}
We will derive the variational inequality from the fact that
\begin{align}
\underset{\theta\to 0}{\lim}\frac{J_1}{\theta}+\underset{\theta\to 0}{\lim}\frac{J_2}{\theta}\geq 0.\label{J1+J2}
\end{align}
Denote by $\left(x_t^\theta,y_t^\theta\right)$ and $\left(x_t^{\left(\theta,\xi\right)},y_t^{\left(\theta,\xi\right)}\right)$ the state trajectories corresponding to $\left(\mu^\theta,\xi^\theta\right)$ and $\left(\mu^\theta,\xi\right)$, respectively. Note that $y_t^{\theta,\xi}=y_t$ since $y_t$ is independent of $\mu_t^{\theta}$.
\begin{lem}\label{prelemma1}
Under assumptions $(A.1)-(A.5)$, we have that for any $p\geq 1$,
\begin{align*}
&\underset{\theta\to 0}{\lim}\E\left|x-x^{\left(\theta,\xi\right)}\right|_T^{*,p}=0,\\
&\underset{\theta\to 0}{\lim}\E\left(\left|x^{\left(\theta,\xi\right)}-x^{\theta}\right|_T^{*,p}+\left|y^{\left(\theta,\xi\right)}-y^{\theta}\right|_T^{*,p}\right)=0.
\end{align*}
\end{lem}
\noindent\textbf{Proof.}
Define $\Delta_t=x_t-x_t^{\left(\theta,\xi\right)}$. Then $\Delta_t$ can be expressed as
\begin{align*}
\Delta_t&=\int_0^t\left(\nu_s\left(\mu_s^\theta\right)-\nu_s\left(\mu_s\right)+x_s^{\left(\theta,\xi\right)}\left(\phi_s\left(\mu_s^\theta\right)-\phi_s\left(\mu_s\right)\right)\right)\d s+\int_0^t\phi_s\left(\mu_s\right)\Delta_s\d s\\
&+\int_0^t\left(\chi_s\left(\mu_s^\theta\right)-\chi_s\left(\mu_s\right)+x_s^{\left(\theta,\xi\right)}\left(\psi_s\left(\mu_s^\theta\right)-\psi_s\left(\mu_s\right)\right)\right)\d B_s+\int_0^t\psi_s\left(\mu_s\right)\Delta_s\d B_s.
\end{align*}
Let
\begin{align*}
z_t=1-\int_0^t\phi_s\left(\mu_s\right)z_s\d s.
\end{align*}
Note that by (A.3) we have that
\begin{align*}
&\E\l|z_\cdot\r|_T^{*,p}<\infty,\\
&\E\l|z^{-1}_\cdot\r|_T^{*,p}<\infty,
\end{align*}
for any $p\geq 1$.
Next we apply Ito's formula to get
\begin{align*}
\Delta_tz_t&=\int_0^t\left(\nu_s\left(\mu_s^\theta\right)-\nu_s\left(\mu_s\right)+x_s^{\left(\theta,\xi\right)}\left(\phi_s\left(\mu_s^\theta\right)-\phi_s\left(\mu_s\right)\right)\right)z_s\d s\\
&+\int_0^t\left(\chi_s\left(\mu_s^\theta\right)-\chi_s\left(\mu_s\right)+x_s^{\left(\theta,\xi\right)}\left(\psi_s\left(\mu_s^\theta\right)-\psi_s\left(\mu_s\right)\right)\right)z_s\d B_s+\int_0^t\psi_s\left(\mu_s\right)\Delta_sz_s\d B_s\\
&=\theta\int_0^t\left(\nu_s\left(q_s\right)-\nu_s\left(\mu_s\right)+x_s^{\left(\theta,\xi\right)}\left(\phi_s\left(q_s\right)-\phi_s\left(\mu_s\right)\right)\right)z_s\d s\\
&+\theta\int_0^t\left(\chi_s\left(q_s\right)-\chi_s\left(\mu_s\right)+x_s^{\left(\theta,\xi\right)}\left(\psi_s\left(q_s\right)-\psi_s\left(\mu_s\right)\right)\right)z_s\d B_s+\int_0^t\psi_s\left(\mu_s\right)\Delta_sz_s\d B_s,
\end{align*}
where we also have used the definition of $\mu^\theta$.
Since $\psi$ is bounded we can apply the Gronwall and Burkholder-Davis-Gundy inequalities to obtain
\begin{align*}
\E|\Delta_\cdot z_\cdot|_T^{*,2p}\leq&~\theta^{2p}K\bigg\{\E\int_0^T\left|\left(\nu_s\left(q_s\right)-\nu_s\left(\mu_s\right)+x_s^{\left(\theta,\xi\right)}\left(\phi_s\left(q_s\right)-\phi_s\left(\mu_s\right)\right)\right)z_s\right|^{2p}\d s\\
&+\E\int_0^T\left|\left(\chi_s\left(q_s\right)-\chi_s\left(\mu_s\right)+x_s^{\left(\theta,\xi\right)}\left(\psi_s\left(q_s\right)-\psi_s\left(\mu_s\right)\right)\right)z_s\right|^{2p}\d s\bigg\}.
\end{align*}
By using the H\"older inequality we can conclude that the integrals on the right hand side are finite, and hence
\begin{align*}
\E\left|x-x^{\left(\theta,\xi\right)}\right|^{*,p}_T=\E\left|z^{-1}_\cdot\Delta_\cdot z_\cdot\right|^{*,p}_T\leq K\left(\E\left|\Delta_\cdot z_\cdot\right|^{*,2p}_T\right)^{1/2}\to 0,
\end{align*}
as $\theta\to 0$.

\medskip \noindent Next, redefine $\Delta$ as $\Delta_t=x_t^{\left(\theta,\xi\right)}-x_t^\theta$. It holds that
\begin{align*}
\Delta_t=&\int_0^t\Delta_s\phi_s\left(\mu_s^\theta\right)\d s+\int_0^t\Delta_s\psi_s\left(\mu_s^\theta\right)\d B_s+\int_0^tG^x_s\d\left(\xi_s^\theta-\xi_s\right)\\
=&\int_0^t\Delta_s\phi_s\left(\mu_s\right)\d s+\int_0^t\Delta_s\psi_s\left(\mu_s\right)\d B_s\\
&+\theta\int_0^t\Delta_s\l(q_s\left(\mu_s\right)-\phi_s\left(\mu_s\right)\r)\d s+\theta\int_0^t\Delta_s\l(\psi_s\left(q_s\right)-\psi_s\left(\mu_s\right)\r)\d B_s\\
&+\theta\int_0^tG^x_s\d\left(\eta_s-\xi_s\right)
\end{align*}
Letting
\begin{align*}
z_t=1-\int_0^t\phi_s\left(\mu_s\right)z_s\d s,
\end{align*}
and applying Ito's formula, yields
\begin{align*}
\Delta_t z_t=&\int_0^t\Delta_s\psi_s\left(\mu_s\right)z_s\d B_s\\
&+\theta\int_0^t\Delta_s\l(\phi_s\left(q_s\right)-\phi_s\left(\mu_s\right)\r)z_s\d s+\theta\int_0^t\Delta_s\l(\psi_s\left(q_s\right)-\psi_s\left(\mu_s\right)\r)z_s\d B_s\\
&+\theta\int_0^tG^x_sz_s\d\left(\eta_s-\xi_s\right).
\end{align*}
Using that $\psi$ is bounded we apply the Burkholder-Davis-Gundy and Gronwall inequalities to obtain
\begin{align*}
\E\left|\Delta_\cdot z_\cdot\right|_T^{*,2p}\leq \theta^{2p}K\bigg\{&\E\int_0^T\left|\Delta_s\l(\phi_s\left(q_s\right)-\phi_s\left(\mu_s\right)\r)z_s\right|^{2p}\d s\\
+&\E\int_0^T\left|\Delta_s\l(\psi_s\left(q_s\right)-\psi_s\left(\mu_s\right)\r)z_s\right|^{2p}\d s\\
+&\E\l|\int_0^Tz_sG^x_s\d\left( \eta_s-\xi_s\right)\right|^{2p}\bigg\}.
\end{align*}
By the H\"older inequality the integral on the right hand side is finite, and hence
\begin{align*}
\E\left|x^{\left(\theta,\xi\right)}-x^\theta\right|_T^{*,p}=\E\left|z^{-1}_\cdot\Delta_\cdot z_\cdot\right|_T^{*,p}\leq K \left(\E\left|\Delta_\cdot z_\cdot\right|_T^{*,2p}\right)^{1/2}\to 0,
\end{align*}
as $\theta\to 0$.
Finally, the fact that $\E\left|y^{\left(\theta,\xi\right)}-y^{\theta}\right|_T^{*,p}\to 0$ as $\theta \to 0$ is a special case of Lemma $3.3$ in \cite{bahlalietal}.
\eop
\begin{lem}
Under assumptions $(A.1)-(A.5)$, we have
\begin{align*}
\underset{\theta\to 0}{\lim}\E\left|\frac{x^\theta-x^{\left(\theta,\xi\right)}}{\theta}-\alpha^x_t\right|_T^{*,2}
+\underset{\theta\to 0}{\lim}\E\left|\frac{y^\theta-y^{\left(\theta,\xi\right)}}{\theta}-\alpha^y_t\right|_T^{*,2}=0,
\end{align*}
where
\begin{align*}
&\alpha^x_t=\int_0^t\phi_s\left(\mu_s\right)\alpha^x_s\d s+\int_0^t\psi_s\left(\mu_s\right)\alpha^x_s\d B_s+\int_0^tG^x_s\d\left(\eta_s-\xi_s\right),\\
&\alpha^y_t=\int_0^tb^y_y\l(s,y_s\r)\alpha^y_s\d s+\int_0^t\sigma_y^y\l(s,y_s\r)\alpha^y_s\d B_s+\int_0^tG^y_s\d\left(\eta_s-\xi_s\right).
\end{align*}
\end{lem}
\noindent\textbf{Proof.}
Let
\begin{align*}
\Delta_t=\frac{x_t^\theta-x_t^{\left(\theta,\xi\right)}}{\theta}-\alpha^x_t.
\end{align*}
Then, we may write $\Delta_t$ as
\begin{align*}
\Delta_t=&\frac{1}{\theta}\bigg\{\int_0^t\phi_s\left(\mu_s^\theta\right)\left(x_s^\theta-x_s^{\left(\theta,\xi\right)}\right)\d s+\int_0^t\psi_s\left(\mu_s^\theta\right)\left(x_s^\theta-x_s^{\left(\theta,\xi\right)}\right)\d B_s\bigg\}\\
&-\int_0^t\phi_s\left(\mu_s^\theta\right)\alpha^x_s\d s-\int_0^t\psi_s\left(\mu_s^\theta\right)\alpha^x_s\d B_s\\
=&\int_0^t\l(\frac{x_s^\theta-x_s^{\left(\theta,\xi\right)}}{\theta}\r)\phi_s\left(\mu_s\right)\d s+\int_0^t\l(\frac{x_s^\theta-x_s^{\left(\theta,\xi\right)}}{\theta}\r)\psi_s\left(\mu_s\right)\d B_s\\
&+\int_0^t\l(x_s^\theta-x_s^{\left(\theta,\xi\right)}\r)\l(\phi_s\left(q_s\right)-\phi_s\left(\mu_s\right)\r)\d s+\int_0^t\l(x_s^\theta-x_s^{\left(\theta,\xi\right)}\r)\l(\psi_s\left(q_s\right)-\psi_s\left(\mu_s\right)\r)\d B_s\\
&-\int_0^t\phi_s\l(\mu_s\r)\alpha_s^x\d s-\int_0^t\psi_s\l(\mu_s\r)\alpha_s^x\d B_s\\
&-\theta\int_0^t\left(\phi_s\left(q_s\right)-\phi_s\left(\mu_s\right)\right)\alpha^x_s\d s-\theta\int_0^t\left(\psi_s\left(q_s\right)-\psi_s\left(\mu_s\right)\right)\alpha^x_s\d B_s\\
=&\int_0^t\Delta_s\phi_s\left(\mu_s\right)\d s+\int_0^t\Delta_s\psi_s\left(\mu_s\right)\d B_s\\
&+\int_0^t\l(x_s^\theta-x_s^{\left(\theta,\xi\right)}\r)\l(\phi_s\left(q_s\right)-\phi_s\left(\mu_s\right)\r)\d s+\int_0^t\l(x_s^\theta-x_s^{\left(\theta,\xi\right)}\r)\l(\psi_s\left(q_s\right)-\psi_s\left(\mu_s\right)\r)\d B_s\\
&-\theta\bigg\{\int_0^t\left(\phi_s\left(q_s\right)-\phi_s\left(\mu_s\right)\right)\alpha^x_s\d s+\int_0^t\left(\psi_s\left(q_s\right)-\psi_s\left(\mu_s\right)\right)\alpha^x_s\d B_s\bigg\}.
\end{align*}
Letting
\begin{align*}
z_t=1-\int_0^t\phi_s\left(\mu_s\right)z_s\d s,
\end{align*}
and applying Ito's formula as well as the Burkholder-Davis-Gundy and H\"older inequalities, we get
\begin{align*}
\E\left|\Delta_\cdot z_\cdot\right|_T^{*,2p}&\leq K\bigg\{\E\int_0^T\left|\Delta_sz_s\right|^{2p}\d s
+\E\int_0^T\l|x_s^\theta-x_s^{\left(\theta,\xi\right)}\r|^{4p}\d s\bigg\}\\
+&\theta\bigg\{\E\int_0^T\left|\left(\phi_s\left(q_s\right)-\phi_s\left(\mu_s\right)\right)\alpha^x_sz_s\right|^{2p}\d s+\E\int_0^T\left|\left(\psi_s\left(q_s\right)-\psi_s\left(\mu_s
\right)\right)\alpha^x_sz_s\right|^{2p}\d s\bigg\}.
\end{align*}
By Lemma \ref{prelemma1} and the integrability of the last term we can apply Gronwall's inequality to conclude that
\begin{align*}
\E\left|\Delta_\cdot z_\cdot\right|_T^{*,2p}\to 0,
\end{align*}
as $\theta \to 0$. Consequently,
\begin{align*}
\E\left|\Delta\right|_T^{*,2}= \E\left|z^{-1}_\cdot \Delta_\cdot z_\cdot\right|_T^{*,2}\leq K \left(\E\left|\Delta_\cdot z_\cdot\right|_T^{*,4}\right)^{1/2} \to 0,
\end{align*}
as $\theta \to 0$.

\medskip\noindent The second assertion is a special case of Lemma 3.4 in \cite{bahlalietal}.
\eop
\begin{lem}
Under assumptions $(A.1)-(A.5)$, we have
\begin{align*}
\underset{\theta\to 0}{\lim}\E\left|\frac{x_\cdot^{\left(\theta,\xi\right)}-x_\cdot}{\theta}-\beta_\cdot\right|_T^{*,2}=0,
\end{align*}
where,
\begin{align*}
\beta_t=&\int_0^t\phi_s\left(\mu_s\right)\beta_s\d s+\int_0^t\psi_s\left(\mu_s\right)\beta_s\d B_s\\
&+\int_0^t\left(b^x\left(s,x_s,\mu_s\right)-b^x\left(s,x_s,q_s\right)\right)\d s+\int_0^t\left(\sigma^x\left(s,x_s,\mu_s\right)-\sigma^x\left(s,x_s,q_s\right)\right)\d B_s.
\end{align*}
\end{lem}
\noindent\textbf{Proof.}
We let
\begin{align*}
\Delta_t=\frac{x_t^{\left(\theta,\xi\right)}-x_t}{\theta}-\beta_t,
\end{align*}
and rewrite $\Delta_t$ as
\begin{align*}
\Delta_t=&\frac{1}{\theta}\int_0^t\left(b^x\left(s,x_s^{\left(\theta,\xi\right)},\mu_s^\theta\right)-b^x\left(s,x_s,\mu_s^\theta\right)\right)\d s\\
&+\frac{1}{\theta}\int_0^t\left(b^x\left(s,x_s,\mu_s^\theta\right)-b^x\left(s,x_s,\mu_s\right)\right)\d s\\
&+\frac{1}{\theta}\int_0^t\left(\sigma^x\left(s,x_s^{\left(\theta,\xi\right)},\mu_s^\theta\right)-\sigma^x\left(s,x_s,\mu_s^\theta\right)\right)\d B_s\\
&+\frac{1}{\theta}\int_0^t\left(\sigma^x\left(s,x_s,\mu_s^\theta\right)-\sigma^x\left(s,x_s,\mu_s\right)\right)\d B_s\\
&-\int_0^t\phi_s\left(\mu_s\right)\beta_s\d s-\int_0^t\psi_s\left(\mu_s\right)\beta_s\d B_s\\
&-\int_0^t\left(b^x\left(s,x_s,\mu_s\right)-b^x\left(s,x_s,q_s\right)\right)\d s-\int_0^t\left(\sigma^x\left(s,x_s,\mu_s\right)-\sigma^x\left(s,x_s,q_s\right)\right)\d B_s\\
=&\frac{1}{\theta}\int_0^t\phi_s\left(\mu_s^\theta\right)\left(x_s^{\left(\theta,\xi\right)}-x_s\right)\d s
+\frac{1}{\theta}\int_0^t\psi_s\left(\mu_s^\theta\right)\left(x_s^{\left(\theta,\xi\right)}-x_s\right)\d B_s\\
&-\int_0^t\phi_s\left(\mu_s\right)\beta_s\d s-\int_0^t\psi_s\left(\mu_s\right)\beta_s\d B_s\\
=&\int_0^t\phi_s\left(\mu_s\right)\Delta_s\d s+\int_0^t\psi_s\left(\mu_s\right)\Delta_s\d B_s\\
&+\int_0^t\left(\phi_s\left(q_s\right)-\phi_s\left(\mu_s\right)\right)\left(x_s^{\left(\theta,\xi\right)}-x_s\right)\d s\\
&+\int_0^t\left(\psi_s\left(q_s\right)-\psi_s\left(\mu_s\right)\right)\left(x_s^{\left(\theta,\xi\right)}-x_s\right)\d B_s.
\end{align*}
Defining
\begin{align*}
z_t=1-\int_0^t\phi_s\left(\mu_s\right)z_s\d s,
\end{align*}
and applying Ito's formula as well as the Burkholder-Davis-Gundy and H\"older inequalities, yields
\begin{align*}
\E\left|\Delta_\cdot z_\cdot\right|_T^{*,2p}\leq&K\left\{\E\int_0^T\left|\Delta_sz_s\right|^{2p}\d s
+\left(\E\int_0^T\left|x_s^{\left(\theta,\xi\right)}-x_s\right|^{4p}\d s\right)^{1/2}\right\}.
\end{align*}
By Gronwall's inequality and Lemma \ref{prelemma1} we have
\begin{align*}
\E\left|\Delta_\cdot z_\cdot\right|_T^{*,2p}\to 0,
\end{align*}
as $\theta\to 0$. Finally, by H\"older's inequality
\begin{align*}
\E\left|\frac{x^\theta-x^{\left(\theta,\xi\right)}}{\theta}-\beta_t\right|_T^{*,2}=\E\left|\Delta\right|_T^{*,2}\to 0,
\end{align*}
as $\theta\to 0$.
\eop
\begin{lem} Under Assumptions $(A.1)-(A.6)$, we have
\begin{align}
\underset{\theta\to 0}{\lim}\frac{J_1}{\theta}&=~\E\left(\alpha^x_Tg_x\left(x_T,y_T\right)+\alpha^y_Tg_y\left(x_T,y_T\right)\right)\notag\\
&+\E\int_0^T\l(\alpha^x_sh_x(s,x_s,y_s,\mu_s)+\alpha^y_sh_y(s,x_s,y_s,\mu_s)\r)\d s\label{J1ineq}\\
&+\E\int_0^Tk_s\d \left(\eta_s-\xi_s\right).\notag
\end{align}
\end{lem}
\noindent\textbf{Proof.}
See \cite{bahlalietal}, Lemma 3.5.
\eop
\begin{lem} Under Assumptions $(A.1)-(A.6)$, we have
\begin{align}
\underset{\theta\to 0}{\lim}\frac{J_2}{\theta}&=~\E\left(g_x\left(x_T,y_T\right)\beta_T\right)\notag\\
&+\E\int_0^Th_x\left(s,x_s,y_s,\mu_s\right)\beta_s\d s \label{J2ineq}\\
&+\E\int_0^T\left(h\left(s,x_s,y_s,q_s\right)-h\left(s,x_s,y_s,\mu_s\right)\right)\d s.\notag
\end{align}
\end{lem}
\noindent\textbf{Proof.}
The proof is similar to that in \cite{peng}, Lemma 2.
\eop


\subsection{Variational inequalities and adjoint equations}
We recall the adjoint processes for the state process (\ref{relaxedstateeq}). These are two pairs of processes $\l(p^x,P^x\r)$, $\l(p^y,P^y\r)$ with values
 in $\mathbb{R}\times \mathbb{R}^{d}$ defined for any control $\l(\mu,\xi\r)\in\mathcal{R}$.
We denote by $f_x$ the derivative with respect to $x$ of the
function $f$, where $f$ stands for either $g$ or $h$. Then $\l(p^x,P^x\r)$, $\l(p^y,P^y\r)$ are given by
\begin{align}
&\left\{ \begin{array}{ll}
\d p^x_t=&-\Big(\phi_t\left(\mu_t\right)p^x_t+\psi_t\left(\mu_t\right)P^x_t+h_x\left(x_t,y_t,\mu_t\right)\Big)\d t+P^x_t\d B_t\\
p^x_T=&g_x\left(x_T,y_T\right) \label{adjoint11}
\end{array} \right.\\
&\left\{ \begin{array}{ll}
\d p^y_t=&-\Big(b^y_y\l(t,y_t\r) p_t+\sigma^y_y\l(t,y_t\r) P^y_t+h_y\left(x_t,y_t,\mu_t\right)\Big)\d t+P^y_t\d B_t\\
p^y_T=&g_y\left(x_T,y_T\right) \label{adjoint12}
\end{array} \right.
\end{align}
Note that the reason for the extra components $P^x,P^y$ is to make it possible to find an adapted solutions to these backward SDEs (see \cite{mayong2} for further discussion). Next, we introduce the Hamiltonian of the system (see e.g. \cite{bensoussan}):
\begin{align*}
&H(t,x,y,\mu,p,P)=-p\Big(\upsilon_t\left(\mu\right)+\phi_t\left(\mu\right)x\Big)-P\Big(\chi_t\left(\mu\right)+\psi_t\left(\mu\right)x\Big)-h\left(t,x,y,\mu\right)
\end{align*}
for $(t,x,y,\mu,p,P)\in[0,T]\times\mathbb{R}\times\mathbb{R}\times \mathcal{P}(U)\times\mathbb{R}\times\mathbb{R}^d$.
\begin{theo} \label{maxprincipleintegral}
(The maximum principle in integral form)~\\ Let $(\mu,\xi)\in \mathcal{R}$ be an optimal relaxed control, i.e
\begin{align*}
J(\mu,\xi)= \underset{(q,\eta)\in \mathcal{R}}{\inf}J(q,\eta),
\end{align*}
and let $\l(x_t,y_t\r)$ be the corresponding trajectory. Then the following inequality holds.
\begin{align}
0\leq& \E\int_0^T\left(k_t+G^x_t\cdot p^x_t+G^y_t\cdot p^y_t\right)d(\eta-\xi)_t\notag\\
&+\E\int_0^T\left(H\left(t,x_t,y_t,\mu_t,p^x_t,P^x_t\right)-H\left(t,x_t,y_t,q_t,p^x_t,P^x_t\right)\right)\d t,\label{intmaxprinciple}
\end{align}%
for all $(q,\eta)\in\mathcal{R}$.
\end{theo}
\noindent\textbf{Proof.}
Define
\begin{align*}
\Phi(x)_t&=1+\int_0^t\phi_s(\mu_s)\Phi(x)_s\d s+\int_0^t\psi_s(\mu_s)\Phi(x)_s\d B_s,\\
\Phi(y)_t&=1+\int_0^tb^y_y\l(s,y_s\r)\Phi(y)_s\d s+\int_0^t\sigma^y_y\l(s,y_s\r)\Phi(y)_s\d B_s\\.
\end{align*}
By Ito's formula $\Phi(x)^{-1},\Phi(y)^{-1}$ are given by
\begin{align*}
\Phi(x)_t^{-1}&=1+\int_0^t\big(\psi_s^T(\mu_s)\psi_s(\mu_s)-\phi_s(\mu_s)\big)\Phi(x)^{-1}_s\d s-\int_0^t\psi_s(\mu_s)\Phi(x)^{-1}_s\d B_s,\\
\Phi(y)_t^{-1}&=1+\int_0^t\big(\sigma^y_y\l(s,y_s\r)^T\sigma^y_y\l(s,y_s\r)-b^y_y\l(s,y_s\r)\big)\Phi(y)^{-1}_s\d s-\int_0^t\sigma^y_y\l(s,y_s\r)\Phi(y)^{-1}_s\d B_s.
\end{align*}
By a simple manipulation we deduce the moment property
\begin{subequations}\label{varineq4}
\begin{align}
&\E|\Phi(x)_\cdot|_T^{*,p}+\E|\Phi(x)^{-1}_\cdot|_T^{*,p}<\infty,\\
&\E|\Phi(y)_\cdot|_T^{*,p}+\E|\Phi(y)^{-1}_\cdot|_T^{*,p}<\infty,
\end{align}
\end{subequations}
for any $p\geq1$. Next, we introduce\begin{align*}
X&=\Phi(x)_Tg_x(x_T,y_T)+\int_0^T\Phi(x)_th_x(t,x_t,y_t,\mu_t)\d t,\\
Y&=\Phi(y)_Tg_y(x_T,y_T)+\int_0^T\Phi(y)_th_y(t,x_t,y_t,\mu_t)\d t,\\
X_t&=\E\big(X|\mathcal{F}_t\big)-\int_0^t\Phi(x)_sh_x(s,x_s,y_s,\mu_s)\d s,\\
Y_t&=\E\big(Y|\mathcal{F}_t\big)-\int_0^t\Phi(y)_sh_y(s,x_s,y_s,\mu_s)\d s.
\end{align*}
Since $g_x$ and $h_x$ are bounded we can use (\ref{varineq4}) to deduce that
\begin{align*}
\E\big|X\big|^p+\E\big|Y\big|^p<\infty,
\end{align*}
for any $p\geq1$. Thus, by the martingale representation theorem (cf.~\cite{karatzasshreve}, Theorem 4.2) there exist $\mathcal{F}_t$-adapted processes $H^x_t,H^y_t$ with the property that, for any $p\geq1$,
\begin{align*}
\E\l(\int_0^T|H^x_t|^2\d t\r)^{p/2}+\E\l(\int_0^T|H^x_t|^2\d t\r)^{p/2}<\infty,
\end{align*}
and such that
\begin{align*}
X_t&=\E(X)+\int_0^tH^x_s\d B_s-\int_0^t\Phi(x)_sh_x(s,x_s,y_s,\mu_s)\d s,\\
Y_t&=\E(Y)+\int_0^tH^y_s\d B^y_s-\int_0^t\Phi(y)_sh_y(s,x_s,y_s,\mu_s)\d s.
\end{align*}
We may now define our adjoint processes $\l(p^x,P^x\r), \l(p^y,P^y\r)$ as
\begin{align*}
p^x_t&=\Phi(x)_t^{-1}X_t,\\
P^x_t&=\Phi(x)_t^{-1}H^x_t-\psi_t(\mu_t)p^x_t,\\
p^y_t&=\Phi(y)_t^{-1}Y_t,\\
P^y_t&=\Phi(y)_t^{-1}H^y_t-\sigma^y_y\l(t,y_t\r) p^y_t,
\end{align*}
noting that $p^x_t$ ($p^y_t$) and $P^x_t$ ($P^y_t$) are $\mathbb{R}$- resp. $\mathbb{R}^d$-valued $\mathcal{F}_t$-adapted processes satisfying
\begin{align*}
\E|p_\cdot|_T^{*,p}+\E\l(\int_0^T|P_\cdot|^2\d t\r)^{p/2}<\infty,
\end{align*}
for any $p\geq1$. Applying Ito's formula on $p^x_t=\Phi(x)_t^{-1}X_t$ and $p^y_t=\Phi(y)_t^{-1}Y_t$ yields
\begin{align}
\d p^x_t&=-\Big(h_x(t,x_t,y_t,\mu_t)+p^x_t\phi_t(\mu_t)+P^x_t\psi_t(\mu_t)\Big)\d t+P^x_t\d B_t,\label{varineq10}\\
\d p^y_t&=-\Big(h_y(t,x_t,y_t,\mu_t)+p^y_tb^y_y\l(t,y_t\r)+P^y_t\sigma^y_y\l(t,y_t\r)\Big)\d t+P^y_t\d B_t.
\end{align}
Using (\ref{varineq10}), and once again by using Ito's formula we can derive
\begin{align*}
\E\Big(p^x_Tx_T\Big)=&\E\int_0^T\Big(p^x_t\big(\nu_t(\mu_t^\theta)+\phi_t(\mu_t^\theta)x_t-\nu_t(\mu_t)-\phi_t(\mu_t)x_t\big)\\
&+P^x_t\big(\chi_t(\mu_t^\theta)+\psi_t(\mu_t^\theta)x_t-\chi_t(\mu_t)-\psi_t(\mu_t)x_t\big)-h_x(t,x_t,y_t,\mu_t)\beta_t\Big)\d t.
\end{align*}
Thus, we can rewrite (\ref{J2ineq}) as
\begin{align*}
\underset{\theta\to 0}{\lim}\frac{J_2}{\theta}=&\E\int_0^T\left(H\left(t,x_t,y_t,q_t,p^x_t,P^x_t\right)-H\left(t,x_t,y_t,\mu_t,p^x_t,P^x_t\right)\right)\d t.
\end{align*}
Similarly, by applying Ito's formula to $p^x_t\alpha^x_t+p^y_t\alpha^y_t$ we may rewrite (\ref{J1ineq}) as
\begin{align*}
\underset{\theta\to 0}{\lim}\frac{J_1}{\theta}=&\E\int_0^T\left(k_s+G^x_s\cdot p^x_s+G^y_s\cdot p^y_s\right)\d \left(\eta_s-\xi_s\right).
\end{align*}
Finally, by combining these equalities with inequality (\ref{J1+J2}) the result follows.\eop
\subsection{Necessary optimality conditions for relaxed controls}
\begin{theo}\label{maxprinciple}
(The relaxed maximum principle)~\\
Let $(\mu
,\xi)$ be an optimal relaxed control and $\l(x_t,y_t\r)$ the
corresponding optimal trajectory, then
\begin{align}
H\left(t,x_t,y_t,\mu_t,p^x_t,P^x_t\right)=\underset{q\in \mathcal{P}(U)}{\sup}H\left(t,x_t,y_t,q_t,p^x_t,P^x_t\right),~\textit{a.e.}, P\textit{-a.s},\label{max1}
\end{align}
\begin{align}
\P\left(\forall t\in[0,T],~ \forall i~ ; \big(k_t^{(i)}+G^{x(i)}_t p^x_t+G^{y(i)}_t p^y_t\big)\geq 0\right)=1,\label{max2}
\end{align}
\begin{align}
\P\left(\sum_{i=1}^d\mathbb{I}_{\{k_t^{(i)}+G_t^{x(i)} p^x_t+G^{x(i)}_t p^x_t>0\}}\d \xi^i=0\right)=1.\label{max3}
\end{align}
\end{theo}

\noindent\textbf{Proof.}
Using the proof in \cite{cadenillashaussmann}, Theorem 4.2, this result can be shown to follow from the maximum principle in integral form (Theorem \ref{maxprincipleintegral}). \eop


\section{A financial application} \label{example}
In this section we give a financial application of the derived stochastic maximum principle. We consider an investor on a market with two investment opportunities, a stock and a portfolio of bonds. The assets in the bond market are non-defaultable bonds, i.e. financial contracts that are bought today and pay a fixed amount at some
 future time, called the maturity time. At each time $t$, the investor is allowed to buy bonds with any time to maturity in $U$, where $U$ is a compact subset of $\mathbb{R}_+$. The dynamics of the bond prices are given by the Heath-Jarrow-Morton model with Musiela parametrization (see \cite{anderssondjehiche}):
\begin{align}
\d p_t(u)&=p_t(u)\l(r_t^0-r_t(u)-v_t(u)\Theta_t\r)\d t+p_t(u)v_t(u)\d B^x_t. \label{bondprice}
\end{align}
Here $B^x_t$ is a Brownian motion, $v_t$ is the integrated volatility process, $\Theta_t$ is the so-called market price of risk, $r_t(u)$ is the forward interest rate and $r_t^0$ is the short rate.

\medskip\noindent Investing in bonds with the price dynamics as above gives the opportunity to, at any time, choose among a continuum of assets, and therefore
we consider measure-valued portfolios. More specifically, we let $q_t\in\boldsymbol{M}$ denote the relative portfolio weights. Then the value of the investment in the bond market can be derived (see \cite{anderssondjehiche}) as an SDE of the form
\begin{align*}
x_t=x_0+\int_0^tx_s\l(r_s^0-v_s\left(q_s
\right)\Theta_s\right)\d s+\int_0^tx_sv_s\left(q_s
\right)\d B^x_s,
\end{align*}
where $x_0$ is the initial capital.

\medskip\noindent The price of a share of a stock is modeled as a geometric Brownian motion and the value of the investment in the stock at time $t$ is then given by
\begin{align*}
y_t=y_0+\int_0^t\lambda y_s \d s+ \int_0^t\rho y_s \d B^y_s,
\end{align*}
where $B^y_t$ is Brownian motion independent of $B^x_t$, $y_0$ is the initial capital and $\lambda$ and $\rho$ are constants.

\medskip\noindent The investor can choose between investing in the bond market, the stock or to consume.  The investors position at time $t$ is $\l(x_t,y_t\r)$ where
 \begin{align*}
 x_t=&x_0+\int_0^tx_s\left(r_s^0-v_s\left(q_s
\right)\Theta_s-c_s\right)\d s+\int_0^tx_sv_s\left(q_s
\right)\d B^x_s+\left(1-K_1\right) \xi^x_t-\xi^y_t,\\
 y_t=&y_0+\int_0^t\lambda y_s\d s+\int_0^t\rho y_s\d B^y_s-\xi^x_t+\left(1-K_2\right)\xi^y_t.
 \end{align*}

The objective of the investor is to choose a consumption/investment strategy consisting of three adapted processes $\l(c_t,q_t,\xi_t\r)$. The consumption process is required to take values in some compact subset of $\mathbb{R}_+$. $\xi_t=\l(\xi_t^x,\xi_t^y\r)$ is nondecreasing and left continuous with right limits. The value of the bonds sold to buy stocks is recorded by $\xi_t^y$ and $\xi_t^x$ records the value of the stocks sold to buy bonds. The constants $0\leq K_1,K_2<1$ account for the proportional transaction costs incurred whenever money is moved between the stock and the bond market.

\medskip\noindent Assuming that our goal is to minimize a cost functional of the form
\begin{align*}
J(q,c,\xi)=\E\l(\int_0^Th(t,x_t,y_t,q_t,c_t)\d t+g(x_T,y_T)\r),
\end{align*}
we get an optimal control problem on the form (\ref{linearcoeff}),(\ref{relaxedstateeq}),(\ref{costfunctional}), with 
\begin{align*}
\mu_t&=q_t\otimes \delta_{c_t},\\
\xi&=\l(\xi^x,\xi^y\r),\\
 \l(\upsilon_t(\cdot),\phi_t(u,c),\chi_t(\cdot),\psi_t(u)\r)&= \l(0,\l(r_t^0-v_t(u)\Theta_t-c\r),0,v_t(u)\r),\\
\l(b^y\l(\cdot,y\r),\sigma^y\l(\cdot,y\r)\r)&=\l(\lambda y, \rho y\r),\\
G^x&=\l(\l(1-K_1\r), -1\r),\\
G^y&=\l(-1, \l(1-K_2\r)\r).
\end{align*}


\subsection{Optimal investment and consumption with transaction costs}
We consider the following cost functional.
\begin{align}
\underset{q,c,\xi}{\sup}\E\l(\int_0^Te^{-\beta t}f\left(c_t
\right)\d t+g\l(x_T,y_T\r)\r).\label{consfunct}
\end{align}
 Assume further that $\l(q_t,c_t,\xi_t\r)$ with corresponding portfolio $\l(x_t,y_t\r)$ is optimal. Using the relaxed maximum principle we may write down the necessary conditions. The adjoint equations becomes
\begin{align*}
&\left\{ \begin{array}{ll}
\d p^x_t=&-\l(\l(r_t^0-v_t(q_t)\Theta_t-c_t\r)p^x_t+v_t\left(q_t\right)P^x_t\r)\d t+P^x_t\d B^x_t\\
p^x_T=&g_x\l(x_T,y_T\r)
\end{array} \right.\\
&\left\{ \begin{array}{ll}
\d p^y_t=&-\l(\lambda p^y_t+\rho P^y_t\r)\d t+P^y_t\d B^y_t\\
p^y_T=&g_y\l(x_T,y_T\r)
\end{array} \right.
\end{align*}
\begin{remark}
The adjoint processes corresponding to standard interest rate model are quite involved and it is difficult to find explicit solutions and therefore difficult to solve this kind of investment/consumption problem explicitly. We refer to \cite{lou} for some worked-out examples of optimal relaxed controls and to \cite{oksendalsulem} for some examples on how to find optimal singular controls.
\end{remark}
\noindent Moreover, the corresponding Hamiltonian is given by
\begin{align}
&H\l(t,x_t,q_t,c_t,p^x_t,P^x_t\r)=\notag\\
&-p^x_t\l(r_t^0-v_t(q_t)\Theta_t-c_t\r)x_t-P^x_tv_t(q_t)x_t-e^{-\beta t}f(c_t). \label{hamilt}
\end{align}
The function $v_t$ corresponds to different interest rate processes.
Choosing the volatility process to be constant, i.e. choosing the Ho-Lee model,
\begin{align*}
\sigma_t(u)=\sigma,
\end{align*}
and consequently
\begin{align*}
v_t(u)=-\sigma u,
\end{align*}
the short rate $r_t^0$ is a Gaussian process. Under a obvious integrability assumption on $\Theta_t$ we then have that (A.1)-(A.3) are fulfilled. Thus by the relaxed maximum principle, the necessary conditions for the optimality of $\l(q_t,c_t\r)$ is that they
maximize (\ref{hamilt}) with
\begin{align*}
v_t(q_t)=-\sigma\int_Uu q_t(\d u).
\end{align*}
Moreover, the optimal time points of transfer between the bond market and the stock is given by the following.
\begin{align*}
\P\left(\forall t\in[0,T]; \big(\l(1-K_1\r) p^x_t - p^y_t\big)\geq 0, \big(-p^x_t +\l(1-K_2\r) p^y_t\big)\geq 0\right)=1,
\end{align*}
\begin{align*}
\P\left(\mathbb{I}_{\{\l(1-K_1\r) p^x_t - p^y_t>0\}}\d\xi^x_t+\mathbb{I}_{\{-p^x_t + \l(1-K_2\r) p^y_t>0\}}\d \xi_t^y=0\right)=1.
\end{align*}
Another choice of volatility process that induces a mean-reverting Gaussian short rate is
\begin{align*}
\sigma_t(u)=\sigma e^{-cu},
\end{align*}
with constants $\sigma$ and $c$. This is the Hull-White model. Similarly, the necessary conditions for optimality is given as above, with
\begin{align*}
v_t(q_t)=\frac{\sigma}{c}\int_U(e^{-cu}-1)q_t(\d u).
\end{align*}

\end{document}